\documentclass [a4paper,12pt]{amsart}
\usepackage [latin1]{inputenc}
\usepackage[all]{xy}
\usepackage{amsmath}
\usepackage{amsfonts}
\usepackage{amssymb}
\usepackage{graphics}
\usepackage[lite]{amsrefs}
\usepackage{amsfonts}
\usepackage{amssymb}
\usepackage{xypic,amsmath,amssymb}
\usepackage{amsmath}
\usepackage[dvips]{graphicx}
\usepackage{graphicx}
\usepackage[lite]{amsrefs}
\newtheorem{teorema}{Theorem}[section]
\newtheorem{lemma}[teorema]{Lemma}
\newtheorem{propos}[teorema]{Proposition}
\newtheorem{corol}[teorema]{Corollary}
\newtheorem{ex}{Example}[section]
\newtheorem{rem}{Remark}[section]
\newtheorem{defin}[teorema]{Definition}
\def\bt{\begin{teorema}}
\def\et{\end{teorema}}
\def\bp{\begin{propos}}
\def\ep{\end{propos}}
\def\bl{\begin{lemma}}
\def\el{\end{lemma}}
\def\bc{\begin{corol}}
\def\ec{\end{corol}}
\def\br{\begin{rem}\rm}
\def\er{\end{rem}}
\def\bex{\begin{ex}\rm}
\def\eex{\end{ex}}
\def\bd{\begin{defin}}
\def\ed{\end{defin}}
\def\demo{\par\noindent{\bf Proof.\ }}
\def\enddemo{\ $\Box$\par\vskip.6truecm}
\def\R{{\mathbb R}}   \def\a {\alpha} \def\b {\beta}
\def\N{{\mathbb N}}     \def\d {\delta} \def\e{\varepsilon}
\def\C{{\mathbb C}}       \def\n{\nu} 
                   \def\r{\varrho} 
                   \def\z{\zeta}
 \def\sbs{\subset} \def\Sbs{\Subset} 
\def\ES{\varnothing}
 \def\oli{\overline} 

  \def\oli{\overline}
\def\G{\Gamma}

\def\oli{\overline}  \def\ES{\varnothing}  \def\IN{\infty} \def\tms{\times} \def\nin{\noindent}

\def\G{\Gamma}

\def\p{\partial}
\def\O{\Omega}
\def\S{\Sigma}

\def\sbs{\subset} \def\sbse{\subseteq}  
  
\def\smi{\setminus} \def\nsbse{\not\subseteqq}

\def\mto{\mapsto}

\def\sps{\supset}
\def\nli{\\}

\def\bit{\begin{itemize}} \def\eit{\end{itemize}}
\def\rar{\rightarrow}
\def\tms{\times}
\def\nid{\not\equiv}

\def\benu{\begin{enumerate}}
\def\eenu{\end{enumerate}}
\def\beqn{\begin{eqnarray*}}
\def\eeqn{\end{eqnarray*}}

\def\bit{\begin{itemize}}
\def\eit{\end{itemize}}
\def\be{\begin{equation}}
\def\ee{\end{equation}}
\def\beq{\begin{eqnarray}}
\def\eeq{\end{eqnarray}}
\def\ba{\begin{array}}
\def\ea{\end{array}}
\def\bt{\begin{teorema}}
\def\et{\end{teorema}}
\def\bp{\begin{propos}}
\def\ep{\end{propos}}
\def\bl{\begin{lemma}}
\def\el{\end{lemma}}
\def\bc{\begin{corol}}
\def\ec{\end{corol}}
\def\br{\begin{rem}\rm}
\def\er{\end{rem}}
\def\bex{\begin{ex}\rm}
\def\eex{\end{ex}}
\def\bd{\begin{defin}}

\newcommand{\hcal}{\mathcal H}

\begin {document}
\baselineskip=.9cm
\title[Some results on evolution]{Some results on evolution}
\maketitle
\tableofcontents
\section{Introduction.}\nin
Let $z_1,\ldots,z_n$ be complex coordinates in $\C^n$, $n\ge 2$. Given a smooth function $\r$ we set 
$$
\vert\partial \r\vert^2=\sum\limits_{1\le\a\le n}\vert \r_\a\vert^2,
$$
$\r_\alpha=\r_{z_\alpha},\r_{\bar\alpha}=\r_{{\bar z}_\alpha}, 1\le \alpha\le n$.

Let $M$ be a smooth hypersurface in $\C^n$ of local equation $\r=0$. For every point $p\in M$ let $HT_p(M)\sbs T_p(M)$ be the complex tangent hyperplane to $M$ at $p$ and $\nu=(\r_{\bar1},\ldots,\r_{\bar n})$ the normal vector to $HT_p(M)$. 

Let $\{E_1,\ldots,E_n\}$ be an ortormal frame with origin at $p$ and such that $\{E_1,\ldots,E_{n-1}\}$ is a frame in $HT_p(M)$ and  $\z_1,\ldots,\z_n$ the complex coordinates determined by $\{E_1,\ldots,E_n\}$. The restriction to $\{\z_n=0\}$ of the Levi form of $\r$ is the {\it intrinsec Levi form} of $M$ at $p$. Its trace is
$$
\hcal(\r)=\vert\p\r\vert^{-1}\sum\limits_{\a,\b=1}^n\left(\d^{\a\bar\b}-\frac{\r_\a\r_{\bar\b}}{\vert\p \r\vert^2}\right)\r_{\a\bar\b}
$$
at $p$.

For $n=2$, $\hcal$ is esentially the Levi operator.

Let $K$ be a compact subset of $\C^n$, $g:\C^n\to\R$ a continuous function which is constant for
$\vert z\vert\gg 0$ and such that $K=\{g=0\}$. Assume that $v\in C^0(\C^n\times \R^+ )$ is a {\em weak solution} of the
parabolic problem  
$$
(\star)\>\>\>\begin{cases}
v_t=\sum\limits_{\a,\b=1}^n\left(\d^{\a\b}-\frac{v_\a v_{\bar\b}}{\vert\p v\vert^2}\right)v_{\a\bar\b}\>\> &{\rm in}\>\> \Omega\times (0,+\IN)\\
v=g\>\> & {\rm on}\>\> \C^n\times\{0\}\\
v={\rm const}\>\>&{\rm for}\>\> t\gg 0.
\end{cases}
$$

Then the family $\{K_t\}_{t\ge0}$ of the subsets $K_t=\{z\in\C^n:v(z,t)=0\}$ (which 
actually depends only on $K$) is called the {\em evolution of $K$ by $\hcal$}.

Evolution of a compact subset $K$ of $\C^2$ was introduced in \cite{ST 7}, \cite{ST 1} where, after proving that the parabolic problem has a unique (weak) solution $u$, it was shown that if $\Omega$ is a bounded pseudoconvex
domain of $\C^2$ with boundary of class $C^3$, the evolution $\{{\overline\Omega}_t\}_{t\ge0}$ of ${\overline\Omega}$ is contained in
${\overline\Omega}$. Conversely, pseudoconcave points "move out by evolution'', i.e. if $\Omega$ is not pseudoconvex then ${\oli\Omega}_t\not\sbse\Omega$ for some $t>0$ (cfr. \cite[Theorem 0.1]{ST 3}). The natural problem of what kind of hull one can recover by evolution was investigated in \cite{ST2}.

In this paper we consider the evolution of a compact subset of $\C^n$ by $\hcal$ {\it with a fixed part} $K^\ast \sbse K$. Precisely, we study the following parabolic problem:
$$
(P)\>\>\>\begin{cases}
v_t=\sum\limits_{\a,\b=1}^n\left(\d^{\a\b}-\frac{v_\a v_{\bar\b}}{\vert\p v\vert^2}\right)v_{\a\bar\b}\>\> &{\rm in}\>\> \Omega\times (0,+\IN)\\
v=g\>\> &{\rm } {\rm on}\>\> {\overline \Omega}\tms\{0\}\\
v(z,t)=g(z)\>\> & {\rm for}\>\> z\in {\rm b}\,\Omega\tms (0,+\IN)
\end{cases}
$$
where $\Omega$ is a bounded strictly pseudoconvex domain in $\C^n$ such that 
$$
K\smi K^\ast\sbse\Omega,\>\>K^\ast\sbse {\rm b}\,\Omega
$$
and $g:\overline\Omega\to\R$ is a continuous function such that $g^{-1}(0)=K$. In Section 2 (see Theorems \ref{EQ12}, \ref {EQ15}) we will prove that
\bit
\item[a)] the problem $(P)$ has a unique (weak) solution $v$ which is bounded and uniformly continuous in ${\overline \Omega}\tms[0,+\IN)$;
\item[b)] if $g$ is a $C^2$ function, the corresponding solution $v$ of $(P)$ is Lipschitz on $\oli\O\tms[0,+\infty)$;
\item[c)] the set
$$
X=\left\{(z,t)\in {\overline \Omega}\tms[0,+\IN):v(z,t)=0\right\}
$$
satisfies
$$
X\cap\left({\overline\O}\tms\{0\}\right)=K\tms\{0\},\>\>X\cap\left({\rm b}\,\Omega\tms[0,+\IN)\right)=K^\ast\tms[0,+\IN)
$$
and it is actually independent of the choice of $g$ and $\Omega$.
\eit

The family $\{E_t(K,K^\ast)\}_{t\ge 0}$ of compact subsets defined by
$$
E_t(K,K^\ast)=\left\{z\in\C^n:(z,t)\in X\> {\rm i.e.}\>\>v(z,t)=0\right\}.
$$
is then said to be the {\em evolution} of $K$ with fixed part $K^\ast$ (by $\hcal$).

Of particular interest in this setting is the case when $K$ is the graph $M$ of a continuous function on the closure $\overline D$ of a bounded domain $D$ in $\C^{n-1}\tms\R$ and $K^\ast={\rm b}\,M$ is the boundary ${\rm b}\,M$ of $M$. 
Generalizing the results of \cite{ST4} for $n=2$ we then prove the following theorem (see Theorem \ref{EQ50}: if $D$ is bounded, strictly pseudoconvex domain i.e. $D\tms i\R$ is a strictly pseudoconvex domain in $\C^n$ then
\begin{enumerate}
\item[d)] $E_t(M,{\rm b}\,M)$ is a graph for all $t\ge 0$ (Theorem \ref{EQ25});  
\item[e)] if ${\rm b}\,M$ is smooth and satisfies the compatibility conditions discovered in \cite{DTZ}, then asymptotically $E_t(M,{\rm b}\,M)$ approaches, in the $C^0$-topology the Levi flat hypersurface with boundary ${\rm b}\,M$ whose existence was proved in \cite{DTZ}.
\end{enumerate}

Let us mention that in the smooth case a parabolic initial value problem related to the flow of a real hypersurface of $\C^n$ by the trace of the Levi form is studied in a nice paper by Huisken and Klingenberg (cfr. \cite{HK}).
\section{Solution of the parabolic problem.}
\subsection{Geometric properties of weak solutions.}
Let $U\sbs\C^n\times (0,+\IN)$ be an open subset. An
upper semicontinuous function $v:U\to [-\IN,+\IN)$ is said to be a {\em (weak) subsolution} of
$$
v_t= \hcal(v)=\sum\limits_{\a,\b=1}^n(\d^{\a\bar\b}-\vert\p v\vert^{-2}v_\a v_{\bar\b})v_{\a\bar\b}.
$$ 
if, for every $(z^0,t^0)$ and a (viscosity) test function 
$\phi$ at $(z^0,t^0)$ (i.e. $\phi$ is smooth near
$(z^0,t^0)$ and $v-\phi$ has a local maximum at $(z^0,t^0)$), one
has  
$$
\phi{_t}(z^0,t^0)\le \hcal(\phi)(z^0,t^0)
$$ 
if $\partial\phi (z^0,t^0)\neq 0$ and
$$
\phi{_t}(z^0,t^0)\le\sum\limits_{\a,\b=1}^n\left(\d^{\a\bar\b}-\eta_\a\eta_{\bar\b}\right)\phi_{\a\bar\b}(z^0,t^0)
$$
for some $\eta\in\C^n$ with $\vert\eta\vert\le 1$, if $\partial\phi(z^0,t^0)=0$. 

A lower semicontinuous function
$v:U\to (-\IN,+\IN]$ is said to be a {\em (weak) supersolution} if, for every
$(z^0,t^0)$ and a test function $\phi$ at
$(z^0,t^0)$  (i.e. $\phi$ is smooth near $(z^0,t^0)$ and
$v-\phi$ has a local minimum at $(z^0,t^0)$), one has  
$$
\phi{_t}(z^0,t^0)\ge \hcal(\phi)(z^0,t^0)
$$
if $\partial\phi (z^0,t^0)\neq 0$ and
$$
\phi{_t}(z^0,t^0)\ge\sum\limits_{\a,\b=1}^n\left(\d^{\a\bar\b}-\eta_\a\eta_{\bar\b}\right)\phi_{\a\bar\b}(z^0,t^0)
$$
for some $\eta\in\C^n$ with $\vert\eta\vert\le 1$, if $\partial\phi(z^0,t^0)=0$.
\br\label{RE1}
Let $A$ be an $n\tms n$ hermitian matrix and $\eta\in\C^n$ with $\vert\eta\vert\le 1$. Then ${\rm Tr}A>{\oli\eta}^{\rm t}A\eta$ provided $A>0$. Conversely, if ${\rm Tr}A>{\oli\eta}^{\rm t}A\eta$ for some $\eta\in\C^n$ with $\vert\eta\vert\le 1$ then $A$ cannot be negative definite. In particular, from the above definition it follows that plurisubharmonic functions are (weak) subsolutions to $v_t=\hcal(v)$.
\er
A {\em (weak) solution} is a continuous function which is both a subsolution and a supersolution.\nli 
One checks that the following properties are true:
\begin{itemize}
\item[1)] maximum (minimum) of a finite number of subsolutions (supersolutions) is a
subsolution (supersolution);
\item[2)] if $W'\subset W\subset\C^n\times (0,+\infty)$, $W$, $W'$ open and 
$v: W\to (-\infty,+\infty)$, $v': W'\to [-\infty,+\infty)$ are subsolutions, such
that for all $\zeta\in bW'\cap W$
$$
\lim\sup_{z\rightarrow\zeta}v'(z)\le v(\zeta)
$$ 
then the function
$$ w(z)=\left\{\begin{array}{ll}\max\;(v(z),v'(z))&\mbox{if $z\in W'$}\\ v(z)&\mbox
{if $z\in W\setminus W'$}
\end{array}\right.
$$ is a subsolution in $W$;
\item[3)] translations of subsolutions (supersolutions) are subsolutions
(supersolutions); i.e. if $\zeta\in \C^n$, $h\in\R$ is positive and 
$v^{\zeta,h}(z,t)=v(z+\zeta,t+h)$ then $v^{\zeta,h}$ is a 
subsolution (supersolution) provided $v$ is;
\item[4)] the limit of a decreasing sequence of subsolutions is a subsolution.
\end{itemize}
\bl\label{EQ1} If $\r: (a,b)\to\R$ is a continuous non decreasing function and   $v$ is a subsolution (or a supersolution) with the range of $v$ in $(a,b)$, then
$\r\circ v$ is a subsolution (or a supersolution, respectively). In particular, if
$v$ is a weak solution, then $\r\circ v$ is a solution.
\el 
\demo There is a sequence of $C^{\infty}$ functions $\r_n: (a,b)\to\R$
such that $\r_n'(t)>0$, $\r_n(t)\searrow\r(t)$, $t\in (a,b)$, therefore it
suffices to prove the lemma for $\r=\r_n$, for then $\r_n\circ v\searrow\r\circ v$
and  $\r\circ v$ will be a subsolution due to 4). Let $\phi$ be a
test function for $\r\circ v$. Then $\chi=\r^{-1}$ is smooth and (strictly) increasing; since $\psi\chi\circ\phi$ is a test function for
$v$ hence we have
$$
\psi{_t}(z^0,t^0)\le \hcal(\psi)(z^0,t^0)
$$ 
if $\partial\psi (z^0,t^0)\neq 0$ and
$$
\psi_t(z^0,t^0)\le \sum\limits_{\a,\b=1}^n\left(\d^{\a\bar\b}-\eta_\a\eta_{\bar\b}\right)\psi_{\a\bar\b}(z^0,t^0)
$$
for some $\eta\in\C^n$ with $\vert\eta\vert\le 1$, if $\partial\psi(z^0,t^0)=0$.

Consider now the case $\partial\psi(z^0,t^0)\neq 0$ and suppose, by a contradiction, that
$$
\phi{_t}(z^0,t^0)>\hcal(\phi)(z^0,t^0).
$$ 
Then
\beqn  
\psi_t(z^0,t^0)&=&\chi'(\phi(z^0,t^0))\phi_t(z^0,t^0)\\
&>&\chi'(\phi(z^0,t^0))\hcal(\phi)(z^0,t^0)=\hcal(\psi)(z^0,t^0)\\
\end{eqnarray*} which is absurd.

As for the case $\partial\psi (z^0,t^0)=0$ it is enough to show the following: let $W\subset\C^n$ 
be open and $\r:W\to\R$ a weak continuous solution of the inequality  
$$
\hcal(\r)(z)\ge-h(z)
$$
where $h:W\to\R^+$ is a continuous positive function. Suppose that $\chi$ is a continuous increasing function $\R\to\R$ with 
$\chi'\in L^\infty(\R)$ and $0\le\chi'\le 1$. Then 
$$
H(\chi\circ\r)(z)\ge-h(z),
$$
in the weak sense. We proceed as follows. Since $\chi$ can be approximated
uniformly on compact subsets of $\R$ by smooth functions with the required
properties, we may assume that $\chi:\R\to\R$, 
$\chi\in C^{\infty}(\R)$, $0<\chi'(s)\le 1$; hence $\chi^{-1}\in
C^{\infty}(\R)$. Let $\psi$ be a smooth test function for
$\hcal(\chi\circ\r)\ge-h$, i.e. 
$$\psi(z)\ge (\chi\circ\r)(z)\>\>{\rm and}\>\>
\psi(z^0)=(\chi\circ\r)(z^0);$$ 
then
$\psi^{\ast}=\chi^{-1}\circ\r$ is a test function too, i.e.
$$\psi^{\star}(z)\ge\r(z), \psi^{\ast}(z^0)=\r(z^0).$$
If
${\partial\psi(z^0)}\neq 0$ we have ${\partial\psi^\ast(z^0)}\neq 0$ and,
by virtue of the hypothesis, $\hcal(\psi^\ast)(z^0)\ge-h(z^0)$, hence

\begin {eqnarray*}
\hcal(\psi)(z)=\hcal(\chi\circ\psi^\ast)(z^0)&=&
\chi'(\psi^{\ast}(z^0))\hcal(\psi^\star)(z^0)\ge\\&\ge&-\chi'(\psi^{\ast}(z^0))h(z^0)>-h(z^0).
\end {eqnarray*}
If $\partial\psi (z^0)=0$, then $\partial\psi^{\ast}(z^0)=0$ and there
is a vector $\eta\in\C^n$, $\vert\eta\vert\le 1$, with
$$
\sum\limits_{\a,\b=1}^n\left(\d^{\a\bar\b}-\eta_\a\eta_{\bar\b}\right)\phi_{\a\bar\b}(z^0,t^0)
\ge-h(z^0).
$$
Now we observe that, since $\psi^{\ast}_\a(z^0)=0$, $1\le\a\le n$
\begin {eqnarray*}
\sum\limits_{\a,\b=1}^n\left(\d^{\a\bar\b}-\eta_\a\eta_{\bar\b}\right)\phi_{\a\bar\b}(z^0,t^0)
&=&
\chi'(\psi^{\ast}(z^0))\phi{_t}(z^0,t^0)\sum\limits_{\a,\b=1}^n\left(\d^{\a\bar\b}-\eta_\a\eta_{\bar\b}\right)\phi_{\a\bar\b}(z^0,t^0)
\ge\\
&&-\chi'(\psi^{\ast}(z^0))h(z^0)\ge
-h(z^0).
\end {eqnarray*}
This ends the proof.
\enddemo
In the sequel we will use the following 
\bp\label{UR}
Let $\{v_\a\}_{\a\in A}$ be a family of weak subsolution of $v_t=\hcal(v)$ and assume that $v=\sup\limits_{\a\in A}v_\a$ is locally bounded from above. Then the upper semicontinuous regularization of $v$
$$
v^\ast(z,t)=\limsup\limits_{(z',t')\to  (z,t)}v(z',t').
$$
is a weak subsolution.
\ep
\demo
We first prove the following: let $B\Sbs W$ be a ball of radius $r$ centered at $w^0=(z^0,t^0)$ and $\phi$ be such that $(v-\phi)(w^0)> (v-\phi)(w)$ for $w\in \oli B\smi{w^0}$. Then there is a sequence $w^\nu\to w^0$ and indices $\a_\nu\in A$ such that for every $\nu$ the function $v_{\a_\nu}-\phi$ has a maximum at $w^\nu$ (relative to $\oli B$).

We may assume that $(v-\phi)(w^0)=0$. For every $\nu\in\N$ such that $1/\nu\le r$ let 
$$
-\d_\nu=\max\left\{(v-\phi)(w):1/\nu\le r\vert w-w^0\vert\le r\right\}.
$$
Since $v-\phi$ has a strict maximum (=0) at $w^0$ (relative to $\oli B$), $-\d_\nu<0$ i.e. $\d_\nu>0$. By definition of regularization
$$
\Big\{(w,s)\in B\tms[-\infty,+\infty):s\le(v^\ast-\phi)(w)\Big\}
$$
is the closure of 
$$
\bigcup_{\a\in A}\big\{(w,s)\in B\tms[-\infty,+\infty):s\le (v_\a-\phi)(w)\big\}.
$$
Thus, for every $\nu$ there is a point $(w^\nu,s^\nu)\in B\tms\R$ and $\a_\nu\in A$ such that 
$$
s^\nu\le (v_{\a_\nu}-\phi)(w^\nu)\le 0,\>\>\vert w^\nu-w^0\vert+s^\nu\le\frac{1}{2}\min(\d_\nu,1/\nu);
$$
in particular
$$
\vert w^\nu-w^0\vert\le\frac{1}{\n},\>\>-\frac{1}{2}\d_\n(v_{\a_\nu}-\phi)(w^\nu)\le 0.
$$
Let now $w^\n$ denote any of the maximum points of $(v_{\a_\nu}-\phi)_{|\oli B}$. Since
\begin {eqnarray*}
(v_{\a_\nu}-\phi)(w^\nu)&\ge&\!\!-\frac{1}{2}\d_\n>-\d_\n\\
&\ge&\!\!\max\Big\{(v-\phi)(w):{\n}^{-1}\le\vert w-w^0\vert\le r\Big\}\\
&>&\!\!\max\Big\{(v_{\a_\n}-\phi)(w):{\n}^{-1}\le\vert w-w^0\vert\le r\Big\}
\end {eqnarray*}
we conclude that $\vert w^\n-w^0\vert\le{\n}^{-1}$ i.e. $w^\n\to w^0$.

In order to prove that $v^\ast$ is a weak subsolution let $\phi\in C^\infty(B)$ and suppose that $v^\ast-\phi$ has a maximum at $w^0=(z^0,t^0)$ with $\p\phi(z^0,t^0)\neq 0$. Let $\phi_\e(w)=\phi(w)+\e\vert w-w^0\vert^2$; then $\p\phi(z^0,t^0)\neq 0$, $\phi_\e$ has a strict maximum at $w^0$ so, in view of what already proved, there are point $w^\n=(z^\n,t^\n)\to w^0=(z^0,t^0)$ and $\a_\n\in A$ such that $(v_{\a_\n}-\phi_\e)$ have maximum at $w^\n$ with $\p\phi_\e(z^\n,t^n)\neq 0$ and
$$
\hcal(\phi_\e)(z^\n,t^\n)\ge\frac{\p\phi_\e}{\p t}(z^\n,t^\n).
$$
Letting $\n\to+\infty$, we get
$$
\hcal(\phi_\e)(z^0,t^0)\ge\frac{\p\phi_\e}{\p t}(z^0,t^0)
$$
and then with $\e\to 0$ 
$$
\hcal(\phi_\e)(z^0,t^0)\ge\frac{\p\phi_\e}{\p t}(z^0,t^0).
$$
The proof if $\p\phi(z^0,t^0)=0$ is similar.
\enddemo 

Finally, in order to prove the independence of the evolution of the pair $(K,K^\ast)$ on $\O$ (see Introduction, $c)$) we discuss a local maximum property of the level sets of a weak solution $v$. 

For an open set $V$ in $\C^n\tms(0,+\IN)$ set
$$
{\mathcal P}_\hcal(V)=\Big\{\psi\in C^2(V):\psi_t\le \hcal(\psi)\Big\}.
$$
\nin
Let $Z$ be a locally closed subset of $V$. We say that $Z$ has {\em local maximum
property }({\em relative to} ${\mathcal P}_\hcal$) if for every open set
$V\Sbs\C^n\tms(0,+\IN)$ such that ${\overline V}\cap Z$ is closed and
${\overline V}$ is compact, and for every $\psi\in {\mathcal P}_\hcal(V')$ where $V'$ is a neighbourhood of
${\overline V}$ it holds:
$$
\max\limits_{{\overline V}\cap Z}\,\psi=\max\limits_{{\rm b}\,V\cap Z}\,\psi.
$$
\bl\label{EQ9}
Let $W\sbse\C^n\tms(0,+\IN)$ be open, $v:W\to\R$ a weak solution of the
$v_t=\hcal(v)$ and $Z=\left\{v=0\right\}$.
Then
\bit
\item[a)] $Z$ has local maximum property; 

\item[b)] for every $c>0$, $Z^c=\big\{(z,t)\in Z:t\le c\big\}$ has local maximum
property. 
\eit
\el
\demo
We first prove the following. Let $v$ be a weak supersolution of $v_t=\hcal(v)$ in $W$. Fix a point $(z^0,t^0)\in W$ and a neighbourhood $V\sbs W$ of $(z^0,t^0)$. Let $\phi\in C^2(V)$ be such that $\phi(z^0,t^0)=v(z^0,t^0)=c$ and
\be\label{EQ33}
\big\{(z,t)\in V:\phi(z,t)>c\Big\}\sbse\Big\{(z,t)\in V:v(z,t)>c\Big\}.
\ee
Then 
$$
\phi{_t}(z^0,t^0)\ge \hcal(\phi)(z^0,t^0)
$$ 
if $\partial\phi (z^0,t^0)\neq 0$ and
$$
\phi{_t}(z^0,t^0)\ge\sum\limits_{\a,\b=1}^n\left(\d^{\a\bar\b}-\eta_\a\eta_{\bar\b}\right)\phi_{\a\bar\b}(z^0,t^0)
$$
for some $\eta\in\C^n$ with $\vert\eta\vert\le 1$, if $\partial\phi(z^0,t^0)=0$.

Observe that, if there exists a non-decreasing continuous function $\r:\R\to\R$ such that $\r(c)=c$ and $\phi(z,t)\le(\r\circ u)(z,t)$ on a neighbourhood of $(z^0,t^0)$, then $\r\circ u$ is still a weak supersolution, so the conclusions concerning $\phi$ are immediate.

In order to construct $\r$ let $N$ be a compact neighbourhood of $(z^0,t^0)$ such that $N\sbs V\sbs W$. Set $\r_1(s)=c$ for $s\le c$. For every $s$ satisfying
$$
c\le s\le s_\infty:=\sup\Big\{v(z,t):(z,t)\in NBig\}  
$$
let 
$$
R_s=\Big\{(z,t):(z,t)\in N:v(z,t)\le s\Big\}.
$$
Since $v$ is lower semicontinuous, the $R_s$'s are compact and $R_s\sbs R_{s'}$ if $s\le s'$. For $c\le s\le s_\infty$ we then define 
$$
\r_1(s)=\max\{\phi(z,t):(z,t)\in R_s\}.
$$
Clearly, $\r_1$ is a non decreasing upper semicontinuous function, $s\mapsto R_s$ being an upper semicontinuous correspondence. Moreover, $\phi(z,t)\le(\r\circ u)(z,t)$. Indeed, assume for a contradiction that $\phi(z,t)>(\r\circ u)(z,t)$. If $\phi(z,t)>(\r\circ u)(z,t)$ this is impossible as $\r_1\ge c$ always. If $\phi(z,t)>c$, by \ref{EQ33}, $v(z,t)>c$. Let $s=v(z,t)$; then $(z,t)\in R_s$ and so $\r_1(s)\ge\phi(z,t),$ i.e.  $\phi(z,t)\le(\r\circ u)(z,t)$. Choose finally a continuous non decreasing function $\r:\R\to\R$ such that $\r\ge\r_1$, $\r(c)=c$. Then $\phi(z,t)\le(\r\circ u)(z,t)$. (Note that $\r$ can be chosen continuous because $\lim\limits_{s\to 0^+}\r_1(s)=c)$.  

Now suppose the claim a) is false, i.e.
$$
\max\limits_{{\overline V}\cap Z}\,\psi>\max\limits_{{\rm b}\,V\cap Z}\,\psi,
$$
for some $\psi\in{\mathcal P}_\hcal(V')$. Then there is $\e>0$ small enough so
that the
function $\psi^\e=\psi-\e t$ still satsfies
$$
\max\limits_{{\overline V}\cap Z}\,\psi^\e>\max\limits_{bV\cap Z}\,\psi^\e,
$$
and, in addition $\psi^\e_t<\hcal(\psi^\e)$ in $\overline V$. Let $(z^0,t^0)$ denote the point where $\psi^\e$ takes maximum
value, say
$M$, relative to ${\overline V}\cap Z$. Clearly $(z^0,t^0)\in V\cap
Z$, and
\begin{eqnarray*}
\Big\{(z,t)\in V:\psi^\e(z,t)>m\Big\}\sbs {V\smi Z}&=&\Big\{(z,t)\in
V:u(z,t)\neq
0\Big\}\\
&=&\Big\{(z,t)\in V:u(z,t)^2>0\Big\}.
\end{eqnarray*}
If we set $\phi=\psi^\e-m$ and $w=u^2$, then $w$ is
still a weak solution of the parabolic problem,
$\phi(z^0,t^0)=w(z^0,t^0)$ and
$$
\left\{(z,t)\in V:\phi (z,t)>0\right\}\sbs \left\{(z,t)\in V:w(z,t)>0\right\}.
$$
Taking into account what proved in the first part we obtain
$$
\psi^\e_t(z^0,t^0)=\phi_t(z^0,t^0)\ge \hcal(\phi)(z^0,t^0)=
\hcal(\psi^\e)(z^0,t^0)
$$
which is a contradiction.

In order to prove b) fix $c>0$ and consider $\psi$ as in definition of
local maximum
property. Let $\r:\R\to\R$ defined by $\r(t)=0$ if $t\le c$, $\r(t)=-(c-t)^3$ if $t>c$
and, for
$N>0$, $(z,t)\in{\overline V}$, let $\psi^N(z,t)=\psi(z,t)+N\r(t)$. Clearly
$\psi^N\in{\mathcal P}_\hcal(V')$ and so, by part a),
$$
\max\limits_{{\overline V}\cap X}\,\psi^N=\max\limits_{{\rm b}\,V\cap X}\,\psi^N.
$$
Observe, however, that
$$
\lim\limits_{N\to+\IN}\psi^N(z,t)=-\IN
$$
if $t>c$ and
$$
\psi(z,t)^N=\psi(z,t)
$$
for $(z,t)\in X^c$, thus
$$
\lim\limits_{N\to+\IN}\max\limits_{{\overline V}\cap X}\,\psi^N=
\max\limits_{{\overline V}\cap X^c}\,\psi.
$$
The same being true for ${\rm b}\,V\cap X^c$, we conclude that
$$
\max\limits_{{\overline V}\cap X^c}\,\psi^N=\max\limits_{{\rm b}\,V\cap X^c}\,\psi^N.
$$
\enddemo
\subsection{Comparison principle. Walsh Lemma in unbounded domains.}
Let us consider the cylinder $Q=\Omega\times(0,h)$ in $\C^n\times{\overline\R}^+$, where $\Omega$ is a bounded domain of $\C^n$
and let 
$$
\Sigma=(\overline\Omega\times \{0\})\cup (b\Omega\times (0,h)).
$$ 
We have the following comparison principle which can be proved arguing as in \cite[Theorem 1.1]{ST1}.  
\bt\label{EQ10} 
Let $v,w\in C^0(\overline Q)$ be respectively a weak subsolution and a
weak supersolution in $Q$. If $v\le w$ on $\oli\S$ then $v\le w$. In particular, $v\le\max\limits_{\oli\S}v$, $w\ge\min\limits_{\oli\S}w.$
\et
We also need the following unbounded version of the Walsh Lemma proved in \cite{ST2}.

Let $W\sbs\R^N$ be a domain with ${\rm b}\,W\neq\ES$ and ${\mathcal F}={\mathcal
F}(\oli W)$ a class of functions $v$ satisfying the following properties:
\begin{itemize}
\item[1)] $v$ is upper semicontinuous in $\oli W$ and $\sup\limits_{\oli
W}\;v=\sup\limits_{{\rm b}\,W}\;v<+\IN$;
\item[2)] for every constant $\a$, $v+\a\in\mathcal F$, if $v\in\mathcal F$;
\item[3)] if $v$ is locally equal to the maximum of finitely many translates of
functions in $\mathcal F$, then $v\in\mathcal F$;
\item[4)] upper semicontinuous regularization of the supremum of a family of
functions in $\mathcal F$ is a member of $\mathcal F$, provided it is uniformly bounded on
$\oli W$;\\
\item[5)] if $W'$ is relatively open in $\oli W$, $v\in{\mathcal F}(\oli W)$,
$v'\in{\mathcal F}({\oli W}')$ and $v'(\z)\le v(\z)$ for $\z\in ({\rm b}\,W')\cap W$ then the
function
$$ w(z)=\left\{\begin{array}{ll}\max\;(v(z),v'(z))&\mbox{if $z\in W'$}\\ v(z)&\mbox
{if $z\in W$}
\end{array}\right.
$$
belongs to ${\mathcal F}(\oli W)$.
\end{itemize}
\bl\label{EQ11}
 Let $g\in C^0({\rm b}\,W)$ be a bounded uniformly continuous function and
$$ 
v(z)=\sup\,\Big\{w(z):v\in{\mathcal F}, w\le g\>\> {\rm on}\>\> {\rm b}\,W\Big\}.
$$ 
Suppose that $v=g$ on ${\rm b}\,W$ and $v$ is uniformly continuous at the points of
${\rm b}\,W$, with the modulus of continuity $\omega(\d)$, $\lim_{\d\rar 0^+}\omega(\d)=0$,
i.e.
$$
\sup\;\Big\{\vert v(z)-g(\z)\vert,\z\in{\rm b}\,W,z\in{\oli W},\vert
z-\z\vert\le\d\Big\}\le\omega(\d). 
$$ Then $v$ is uniformly continuous on $\oli W$ with the same $\omega(\d)$ as its
modulus of continuity.
\el
\subsection{Existence of solutions and evolution.}
We are in position to prove the following existence theorem:
\bt\label{EQ12}
Let $\Omega$ be a bounded strictly pseudoconvex domain in $\C^n$, $g:\overline\Omega\to\R$ a continuous function. Then the problem $(P)$ has a unique weak solution $v$ which is bounded and uniformly continuous in ${\overline
\Omega}\tms[0,+\IN)$.
\et
\demo
Unicity is a consequence of the comparison principle. Existence will be proved by Perron method.

Let $W=\Omega\tms[0,+\IN)$ and ${\mathcal F}={\mathcal F}_{\rm g}$ be the class  of all functions 
$w:{\overline W}\to [-\IN,+\IN)$ with the following properties:\vskip0.4truecm
1) $w$ is upper semicontinuous in $\overline W$ and is a subsolution in
$W$;\vskip0.2truecm

2) $w\le \max\limits_{\oli W}g$;\vskip0.2truecm

3) $w\le g$ on ${\rm b}\,W$\vskip0.4truecm\nin 
Let $v:\oli W\to\R$ be the function
$$ 
(z,t)\longrightarrow\sup\,\Big\{w(z,t):w\in{\mathcal F}\Big\}
$$ 
and $v^\ast$ its upper semicontinuous regularization: $v^\ast$ is a subsolution (cfr. Proposition \ref{UR}).

We want to prove that $v=v^\ast$ and $v$ is actually the solution of the problem $(P)$. 

The proof is divided in several steps.

A) $v=g$, for $(\z,t)\in {\rm b}\,\Omega\times [0,+\IN)$. Furthermore, $v$ is uniformly continuous at the
points of ${\rm b}\,\Omega\times [0,+\IN)$ in the following sense: given $\e>0$ there exists
$\d>0$ such that $\vert v(z,t)-g\vert\le\e$ if dist\,$\left((z,t),{\rm b}\,\Omega\times
[0,+\IN)\right)<\d$.

Let $\e>0$ be fixed and $g_1\in C^2(\C^n)$ such that $\vert g_1-g\vert<\e$ on $\oli\O$. Since $\Omega$ is strictly pseudoconvex there is a strictly plurisubharmonic function $\r$ on a neighbourhood $U$ of 
$\overline\Omega$, such that
$\r=0$ on ${\rm b}\,\Omega$ and $\O=\{\r<0\}$. For $m>0$ big enough the time-independent function
$$ v^\e_m(z)=m\r(z)+g_1-\e,
$$
is strongly plurisubharmonic in $U$, therefore a subsolution of $v_t=\hcal(v)$ and 
$$
v^\e_m(z)=n\r(z)+g_1-\e\le n\r(z)+g\le g(z)
$$
for $z\in\overline W$.. 

Thus  $v^\e_m\in{\mathcal F}$ and consequently
$$  
m\r(z)+g-\e \le v^\e_m(z,t)\le v(z,t)\le g.
$$
It follows 
$$  
\vert v(z,t)-g\vert \le m\vert\r\vert+\e.
$$
It is evident now that, for a fixed $\e$, there is $\d>0$ such that the statement A) holds true. 

B) For all $a\in \O$
\be\label{EQ13}
\lim\limits_{(z,t)\to(a,0)}v(z,t)=\lim\limits_{(z,t)\to(a,0)}v^\ast(z,t)=g(a)
\ee
In order to prove this we fix $\e>0$ and smooth functions $\phi$, $\psi$ on $\C^n$ in such a way to have
$$
g(z)-\e<\phi(z)<g(z)<\psi(z)<g(z)+\e.
$$
Let $c$ be a constant such that
$$
\left\vert\sum\limits_{\a,\b=1}^n\phi_{\a\bar\b}(z)\xi^\a\oli{\xi^\b}\right\vert<c\vert\xi\vert^2, \left\vert\sum\limits_{\a,\b=1}^n\psi_{\a\bar\b}(z)\xi^\a\oli{\xi^\b}\right\vert<c\vert\xi\vert^2
$$
for all $z\in\O$ and $\xi\in\C^n$. Then
$$
v_+(z,t)=\psi(z)+ct,\>\>v_-(z,t)=\phi(z,t)-ct
$$
are respectively a regular supersolution and a regular subsolution in $W$; moreover, $v_-\in\mathcal F$ and $v_{|{\rm b}\,W}\le v_{+|{\rm b}\,W}$. In view of the comparison principle for $v_+$ and $v_-$, we deduce that
$$
v_-\le v\le v^\ast\le v_+
$$
in $\oli W$ and consequently, since $v_-$ and $v_+$ are continuous, that
$$
g(a)-\e\le v_-(a,0)\le\liminf_{(z,t)\to(a,0)}v(z,t)\le\limsup\limits_{(z,t)\to(a,0)}v(z,t)\le u_+(a,0)\le g(a)+\e.
$$
for all $a\in\oli\O$.

\ref{EQ13} follows $\e$ being arbitrary.

 A), B) imply that $v^\ast\in\mathcal F$ therefore, by definition of $v$, we have $v^\ast=v$. In particular, $v=v^\ast$ is a subsolution which is continuous at every point of ${\rm b}\,W$. Thus all the hypothesis of the Walsh Lemma (see \ref{EQ11}) are satisfied hence $v$ is continuous in $\oli W$.

Finally $v$ is a weak solution in $W$. For if not there is $(z^0,t^0)\in W$ and $\phi\in C^\infty(W)$ such that $v-\phi$ has a strict local minimum (=0) at $(z^0,t^0)$ and
$$
\phi{_t}(z^0,t^0)< \hcal(\phi)(z^0,t^0)
$$
if $\partial\phi (z^0,t^0)\neq 0$ and
\be\label{EQ14}
\phi{_t}(z^0,t^0)<\sum\limits_{\a,\b=1}^n\left(\d^{\a\bar\b}-\eta_\a\eta_{\bar\b}\right)\phi_{\a\bar\b}(z^0,t^0)
\ee
for some $\eta\in\C^n$ with $\vert\eta\vert\le 1$, if $\partial\phi(z^0,t^0)=0$. Observe that $v(z^0,t^0)<\max\limits_{\oli W}g$, otherwise, by definition of $v$, $(z^0,t^0)$ would be a maximum point for $v$ hence for $\phi$ and this contadicts \ref{EQ14} (see Remark \ref{RE1}). Thus, we can find $\e>0$ small enough such that $\phi+\e$ is a subsolution on a neighbourhood $U$ of $(z^0,t^0)$, $\phi+\e<\max\limits_{\oli W}g$ and 
$$
\ES\neq V=\Big\{(z,t)\in U:(\phi+\e-v)(z,t)>0\Big\}\Sbs U.
$$
It is now clear that
$$ 
\tilde v(z,t)=\left\{\begin{array}{ll}\max\;(v(z,t),\phi(z,t)+\e)&\mbox{if
$(z,t)\in U$}\\ v(z,t)&\mbox {if $(z,t)\in {\overline W}\smi \oli V$}
\end{array}\right. 
$$ 
is a subsolution, $\tilde v\in {\mathcal F}$ and $v<\tilde v$ near $(z^0,t^0)$: contradiction. 

Theorem \ref{EQ12} is completely proved.
\enddemo
\br\label{RE2}
The strict pseudoconvexity condition can be relaxed. In particular the following condition suffices: for all $\zeta\in{\rm b\O}$ there is a ball $B$ centered at $\zeta$ and a strictly plurisubharmonic function $\phi:B\to\R$ such that $\phi(\z)=0$ and $\phi<0$ on $B\cap\O$. 
\er
\br\label{RE3}
Using the method employed in \cite{SIM} it can be proved that if the boundary value $g$ is $C^2({\rm b}\,\O)$ the solution of the problem $(P)$ is Lipschitz in $\oli\O.$
\er
\bt\label{EQ15}
Let $(K^\ast,K)$ be a pair of compact sets in $\C^n$ such that $K^\ast\sbs K$, $K\nid\ES$ and $\Omega$ a bounded strictly pseudoconvex domain
such that $K\smi K^\ast\sbse\Omega$, $K^\ast\sbse {\rm b}\,\Omega$. Assume that $K=g^{-1}(0)=K$ with $g:\overline\Omega\to\R$ and let $v$ be the solution of the parabolic problem $(P)$. Then the set 
$$
X=\left\{(z,t)\in {\overline \Omega}\tms[0,+\IN):v(z,t)=0\right\}
$$
is independent of the choice of $g$ and $\Omega$. Moreover
\bit
\item[i)]
$
X\cap\left({\overline\O}\tms\{0\}\right)=K\tms\{0\},
$
\item[ii)]
$
X\cap\left({\rm b}\,\Omega\tms[0,+\IN)\right)=K^\ast\tms[0,+\IN).
$
\eit
\et
\demo
The independence of the zero set $\{u=0\}$ of the choice of $g$ satisfying $g^{-1}(0)=K$ is essentially the argument of Evans and Spruck in \cite{ES1} (cfr. also \cite{ST2}).

It remains to show independence of $X=u^{-1}(0)$ of the choice of $\O$ satisfying the conditions of Theorem \ref{EQ15} . 

Suppose $\O_1$, $\O_2$ are such domains and
$\O_0=\O_1\cap\O_2$. Then $\O_0$ satisfies condition $(C)$ of Remark \ref{RE2} and also
$K\smi K^\ast\sbse\O_0$, $K^\ast\sbse {\rm b}\,\O_0$. For each of these sets we have unique (independent of respective $u$)
"evolution hypersurface" i.e $X_j$, where $j=0, 1, 2$, $X_j\sbse
{\overline\O_j}\tms[0,+\IN)$ and
\bit
\item[i)]
$
X_j\cap\left({\overline\O}\tms\{0\}\right)=K\tms\{0\},
$
\item[ii)]
$
X_j\cap\left({\rm b}\,\Omega\tms[0,+\IN)\right)=K^\ast\tms[0,+\IN).
$
\eit
We will show that $X_1=X_0$ and this will imply that $X_1=X_2$, as required.

Let $g$, $v$ be as in Theorem \ref{EQ15}, for the domain $\O_1$, so that
$X_1=u^{-1}(0)$. Let now $g_0=g_{|{\overline\O_0}}$ and $u_0\in
C^0\left({\overline\O_0}\tms[0,+\IN)\right)$ be
the corresponding solutions of the parabolic problem so that
$X_0=u_0^{-1}(0)$.\nli
The following is true:
\bit
\item[i)] $X_1\sbse {\overline\O_0}\tms[0,+\IN)$;
\item[ii)] $X_1\sbse \left(\O_0\cup K^\ast\right)\tms[0,+\IN)$.
\eit
Since $\O_0$ is the intersection of two strictly pseudoconvex domains $\O_1$,
$\O_2$, there is a neighbourhood $N$ of ${\overline\O_0}$ and a continuous
plurisubharmonic function $\phi:N\to\R$ such that
${\overline\O_0}=\{\phi\le 0\}$. Suppose
$X_1\nsbse{\overline\O_0}$, then there exists $c>0$ such that
$X^c\nsbse{\overline\O_0}$ but
$X^c_1\sbs N$. (Observe that $c\mto X^c$ is an upper semicontinuous
correspondence and
$X^\circ_1=K\sbs N$.)

Let $\widetilde {\phi}(z,t)=\phi(z)$ and define
$$
M=\max\limits_{X^c}\,\widetilde {\phi},\>\>
F=\left\{(z,t)\in X^c:\widetilde
{\phi}(z,t)=M\right\}.
$$
Then $M>0$, $F$ is compact and $F\cap\left(K^\ast\tms\{0\}\right)=\ES$.
Choose $V$, a neighbourhood of $F$ such that $\overline V$ is compact,
$\overline V\sbs
N\smi{\overline\O_0}\tms(0,+\IN)$. Then
$$
M=\max\limits_{X^c_1\cap{\overline V}}\,\tilde
{\phi}>\max\limits_{X^c_1\cap {\rm b}\,V}\,\tilde
{\phi}
$$
which contradicts the local maximum property (b) of Lemma \ref{EQ9} since, clearly, $u\in {\mathcal P}_\hcal$. . This proves i).

As for ii) suppose $(z^0,t^0)\in X_1\cap\left({\rm b}\,\O_0\smi
K^\ast\right)\tms[0,+\IN).$
Then $z^\ast\in {\rm b}\,\O_1$ or $z^\ast\in {\rm b}\,\O_2$. In either case there is a $C^2$
strictly plurisubharmonic function $v=v(z)$ in a neighbourhood of $z^\ast$
such that
$v(z^\ast)=0$, $v(z)<0$ for $z\in
B(z^\ast,r)\cap({\overline\O_0}\smi\{z^\ast\})$. Since
$v$ is strictly plurisubharmonic, there is an $\e>0$, small enough so that
the function
$\psi^\e(z,t)=v(z)-\e(t-t^\ast)^2$ is of the class ${\mathcal P}_\hcal$ in
$V=B\tms(t^0-r,t^\ast+r)$. Observe now that
$\psi^\ast(z^0,t^0)=0$ while
$\psi^\e(z,t)<0$ for $(z,t)\in X_1\cap V\smi\{(z^0,t^0)\}$. This
contradicts
again the local maximum property (a) of Lemma \ref{EQ9}. Thus
$X_1\cap\left({\rm b}\,\O_0\smi K^\ast\right)\tms[0,+\IN)=\ES.$, whence ii).

We can  now show that $X_0=X_1$. Fix $c>0$ and let
$$
W^c=\O_0\tms(0,c),\>\>
\Sigma^c=\left({\overline\O_0}\tms\{0\}\right)\cup\left(b\O_0\tms[0,c]\right).
$$
Let $U^c=u_{|{\overline W^c}}$. Then $u_0$, $U^c$ are continuous weak
solutions in
$\overline W$. By i), ii)
$$
u^{-1}(0)\cap\Sigma^c=\left(U^c\right)^{-1}(0)\cap\Sigma^c.
$$
Hence, similarly as in [ES1] there are continuous increasing functions
$\chi_1,\chi_2:\R\to\R$, with $\chi_j(0)=0$, $j=1, 2$, such that
$$
\chi_1\circ u_0\le U^c\le\chi_2\circ u_0
$$
on $\Sigma^c$.\nli
Since $\chi_j\circ u_0$, $j=1,2$, are weak solutions the comparison principle implies that
$$
\chi_1\circ u_0\le U^c\le\chi_2\circ u_0
$$
in $W$ and so
$$
X^c_0=(u_0)^{-1}(0)=(U^c)^{-1}(0)=X^c_1,
$$
for every $c>0$. Thus $X_0=X_1$.
\enddemo

In light of this theorem we define
$$
E_t(K,K^\ast)=\Big\{z\in\C^n:(z,t)\in X\> {\rm i.e.}\>\>u(z,t)=0\Big\}.
$$
The family $\{E_t(K,K^\ast)\}_{t\ge 0}$ is said to be the {\em evolution} of $K\bmod\,K^\ast$ ({\em by $\hcal$}).

The semigroup property
\be\label{SGR}
E_{t+t'}(K,K^\ast)=E_t\left(E_{t'}(K,K^\ast),K^\ast\right)
\ee
holds true as well as for the standard evolution (i.e. when $K^\ast=\ES$).  
\subsection{Some geometric properties}
\bt\label{EQ16}
Let $\Omega$ be a bounded, strictly pseudoconvex domain of $\C^n$,
$K\sbs{\overline\Omega}$, $K^\ast\sbs {\rm b}\,\Omega$ compact sets such that:
$K^\ast\sbs K$,
$K\smi K^\ast\sbs \Omega$ and separates $\Omega$. Let
$\big\{E_t(K,K^\ast)\big\}_{t\ge 0}$ be
the evolution of $K\,\bmod\,K^\ast$. Then, for every $t$ the subset
$E_t(K,K^\ast)\smi K^\ast$ separates $\Omega$.
\et
\demo
Choose $g\in C^0(\overline\Omega)$ such that
$g^{-1}(0)=K$; $\Omega\smi K=\{g>0\}\cup\{g<0\}$ and we choose
$\z_1,\z_2$ such that $g(\z_1)>0$, $g(\z_2)<0$. Let $u$ be the weak
solution of $(P)$.
Then $E_t(K,K^\ast)=\{u(\cdot,t)=0\}$ and $\Omega\smi
\{E_t(K,K^\ast)\}_{t\ge 0}$ is a union
$\{u(\cdot,t)>0\}\cup\{u(\cdot,t)<0\}$ of
nonempty subsets.
\enddemo
\bp\label{EQ17}
In the context of the previous theorem
$$
\limsup\limits_{t\to+\IN}\, E_t(K,K^\ast)=K^\IN
$$
where $K^\IN\smi K^\ast$ is pseudoconcave i.e. has local maximum property
with respect
to the functions $|P|$, $P\in\C[z_1,z_2,\ldots,z_n]$. Furthermore, $K^\IN\smi K^\ast$
separates
$\O$. 
\ep
\demo
First of all we point out the following fact whose proof is a straightforward consequence of the definition. Let $\left\{X_t\right\}_{t\in T}$, where $T$ is a (direct) partially ordered set, be a family of relatively closed subsets of an open subset $W$ of $\C^n\tms(0,+\IN)$. Assume all $X_t$ have local maximum property relative to ${\mathcal P}_\hcal$. Then
$$
\limsup\limits_{t\to+\IN}\,X_t=\bigcap_{t^\circ}\overline
{\bigcup\limits_{t\ge t^\circ}\,X_t}
$$
has local maximum property relative to ${\mathcal P}_\hcal$ provided it is nonempty.

In order to prove that $K^\IN\smi K^\ast$ is pseudoconcave let $W=\O\tms (0,+\IN)$ and $u$ be the solution of the parabolic problem (P). We know that $v$ is uniformly continuous in $\overline W$.

Let
$$
X=\Big\{(z,t)\in\O\tms(0,+\IN):v(z,t)=0\Big\}
$$
and
$$
X^h=\Big\{(z,t)\in\O\tms(0,+\IN):v^h(z,t)=0\Big\}
$$
where $v^h(z,t)=v(z,t+h)$, $h>0$.

Since the equation $v_t=\hcal(v)$ is invariant with respect to
time shift $t\mto t+h$, $h\ge 0$, we obtain that $\Big\{X^h\cap
W\Big\}_{h>0}$ is a
family of sets with local maximum property relative to ${\mathcal P}_\hcal$
defined above. Let
$$
X^\IN=\limsup\limits_{h\to+\IN}\,X^h.
$$
By what observed at the beginning, $X^\IN\cap W$ has local maximum property relative to ${\mathcal P}_\hcal$ provide $X^\IN\cap W\neq\ES$. On the other hand, from
$$
X^h\cap\left(\C^n\tms\{t\}\right)=E_{t+h}(K,K^\ast)\tms\{t\},
$$
and
$$
\limsup\limits_{h\to+\IN}\,E_{t+h}(K,K^\ast)=K^\IN,
$$
for each $t>0$ we deduce that $X^\IN=K^\IN\tms(0,+\IN)$ and so the set $\left(K^\IN\smi K^\ast\right)\tms(0,+\IN)$ has local maximum property
relative to the class of subsolutions ${\mathcal P}_\hcal$.

Suppose now that $K^\IN\smi K^\ast$ is not a local maximum set relative to the functions $|P|$, $P\in\C[z_1,z_2,\ldots,z_n]$. Then, by \cite{S} there are a point $z^0\in K^\IN\smi K^\ast\sbs\O$ and a strictly
plurisubharmonic function $\r\in C^2\left(B(z^0,r)\right)$, $r>0$, such that
$\r(z^0)=0$ and $\r(z)<0$ for $z\in
K^\IN\cap\left(B(z^0,r)\smi\{z^0\}\right)$.
Choose a small $\e>0$ such that the function
$\psi(z,t)=\r(z)-\e(t-t^0)^2$ satisfies
${\hcal}(\psi)-\psi_t>0$ in $B(z^0,r)\tms(t^0-r,t^0+r)$, i.e.
$\psi\in{\mathcal P}_{\hcal}$ in a neighbourhood of $(z^0,t^0)$. Owing to the
properties of $\r$,
$$
\Large{\psi_{|\left(K^\IN\smi
K^\ast\right)\tms(0,+\IN)}=\psi_{|X^\IN\cap W}}
$$
$$
{\psi_{|\left(K^\IN\smi
K^\ast\right)\tms(0,+\IN)}=\psi_{|X^\IN\cap W}}
$$
has strict local maximum at $(z^0,t^0)$: contradiction.
\enddemo
We will need the following general fact.
\bp\label{EQ19}
Let $\Omega\sbs\C^n$ be a bounded strictly pseudoconvex domain and $K,K'$ disjoint compact subsets of $\O$. Let 
$K\cap{\rm b}\,\O=K^\ast$, $K'\cap {\rm b}\,\O={K'}^\ast$. Then
$$
E_t(K,K^\ast)\cap E_t(K',{K'}^\ast)=\ES
$$
for every $t>0$.  
\ep
\demo
Take a continuous function $g:{\overline\Omega}\to\R$ such that
$g^{-1}(0)=K$, $g^{-1}(1)=K'$ and solve the problem
$$
\left\{\ba{ll} v_t&=\hcal(v) \>\>{\rm in}\>\> \Omega\times (0,+\IN)\\
v&=g\>\> {\rm on}\>\>\left(
{\overline \Omega}\tms\{0\}\right)\cup \left({\rm b}\,\Omega\tms (0,+\IN)\right).
\ea\right.
$$
Then 
$$
E_t(K,K^\ast)=\{v(\cdot,t)=0\},\>\>E_t(K',{K'}^\ast)=\{v(\cdot,t)=1\}
$$
and consequently the subsets $E_t(K,K^\ast)$, $E_t(K',{K'}^\ast)$ are disjoint for every
$t>0$.
\enddemo
\br
We do not know if the same is true if we have two different strictly
pseudoconvex domains $\O$, $\O'$ with $K\sbse\O$, $K'\sbse\O'$.
\er
\section {Evolution of graphs}\nin
From now on we assume that $K$ is the graph $\Gamma$ of a continuous function $u:{\overline D}\to\R$ where $D$ is a bounded domain of $\C^{(n-1)}\tms\R$ and $K^\ast={\rm b}\,\Gamma$. 

We have the following theorem:  
\bt\label{EQ25} 
If $D\tms i\R$ is strictly pseudoconvex then $E_t(\Gamma,{\rm b}\,\Gamma)$ is a graph for every $t\ge 0$.
\et
\demo In our situation $K=\Gamma$ and $D\tms i\R$ is a
strictly pseudoconvex domain in $\C^n$. Set $z'=(z_1,\dots,z_{n-1})$ and consider translations $T_h:\C^n\to\C^n$ of the form
$(z',z_)\mto (z',z_n+ih)$, $h\in\R$.

For fixed $h>0$, consider a bounded strictly pseudoconvex domain $\O$ and
a large
enough number $M$ such that
\begin{eqnarray*}
{\overline D}\tms i\R\sps {\overline\O}&\sps&{\overline D\tms[-iM,iM]}\\
&\sps&\Gamma\cup T_h(\Gamma).
\end{eqnarray*}
Then we can consider the evolutions of $\Gamma$ and $T_h(\Gamma)$
($\bmod\,{\rm b}\,\Gamma$)
with such $\overline\O$ and they must be disjoint in view of Proposition \ref{EQ19}. (The evolution is
independent of the
specific choice of such $\O$.)
\enddemo
The operator $\hcal$ does not depend upon the equation of a surface. In particular, if $x_1, y_1,\dots,x_n,y_n$ are real coordinates with $z_\a=x_\a+iy_\a$, $1\le\a\le n$, for a graph of a smooth function $y_n=u(x_1,\ldots,x_n,y_1,\dots,y_{n-1})$  one has $\hcal(y_n-u)=\hcal_0(u)$ where $\hcal_0$ is a quasilinear degenerate elliptic operator in the real coordinates. 

If $n=2$ $\hcal_0$ is the Levi operator for graphs  (cfr. \cite{ST7})
\begin{eqnarray*}
{\hcal}_\circ(u)&=&\frac{1}{4}(1+\vert
Du\vert^2)^{-1}\big\{(1+u^2_3)(u_{11}+u_{22})+(u^2_{1}+u^2_{2})u_{33}\\
&+&2(u_2-u_1u_3)u_{13}-2(u_1+u_2u_3)u_{23}\big\}
\end{eqnarray*}
($u_j={\partial u}/{\partial x_j}$, $u_{ij}={\partial^2 u}/\partial
x_i \partial x_j$).
\bl\label{EQ20}
Let $u$ be continuous in a domain $D\sbse\C^{n-1}\tms\R$. Then $y_n-u$ is a
weak solution of $v_t={\hcal}(v)$ in $D\tms i\R\tms(0,+\infty)$ if and only if $u$
is a weak solution of $u_t={\hcal}_0 (u)$ in $D\tms\R$.
\el
\demo
Set $x=(x_1,\ldots,x_{n})$, $y'=(y_1,\ldots,_{n-1})$. If $v=y_{n}-u(x,y')$ is a weak solution of $v_t={\hcal}(v)$ in $D\tms i\R\tms(0,+\infty)$
then is
immediately seen that $u$ is a weak solution of $u_t={\hcal}_0(u)$ in
$D\tms(0,+\infty)$.

Conversely, let us suppose that $u$ is a weak
solution of $u_t={\hcal}_0(u)$ and let $\phi=\phi(x,y,t)$ be smooth and
such that
$y_n-u-\phi$ has a local maximum at $(\bar x,\bar y',\bar t)$. We may assume that
$(\bar x,\bar y',\bar t)=(0,0,0)$ and
$u(0,0,0)=\phi(0,0,0)=0$. Since, locally at $(0,0,0)$, $y_n-u\le\phi$ we have
$\phi_{y_n}(0,0,0)=1$. In particular, $\phi=0$ is a (local) graph $y_n=f(x,y',t)$ and
$\phi=\lambda(y_n-h)$ with $\lambda$ smooth and $\lambda(0,0,0)=1$.
Moreover, since ${\hcal}$ is invariant with respect to unitary transformations of
$\C^n$, we may also assume that ${\rm d}_{x,y'}f(0,0,0)=0$. In this situation we have
$$
-u(x,y',t)\le-\lambda(x,0,t)f(x,y',t)
$$
and
$$ {\hcal}_0(-\lambda
f)(0)=-\frac{1}{4}\sum\limits_{j=1}^{n_1}\big[f_{x_jx_j}(0,0,0)+f_{y_jy_j}(0,0,0)\big]=\hcal_0(-f)(0,0,0).
$$
Furthermore
$$
\phi_t(0,0,0)=f_t(0,0,0),\>\> {\hcal}(\phi)(0,0,0)={\hcal}_0(-f)(0,0,0).
$$
Since $-u$ is a weak solution of $w_t={\hcal}_0(w)$
$$
-f_t(0,0,0)\le{\hcal}_(-f)(0,0,0).
$$
From this, in view of the above identities, we obtain
$$
\phi_t(0,0,0)=-f_t(0,0,0)=\le-{\hcal}_0(-f)(0,0,0)= \hcal_0(-f)(0,0,0)=\hcal(\phi)(0,0,0).
$$
This proves that $y_n-u$ is a weak subsolution.

Similarly we prove that $y_n-u$ is a  weak supersolution.

Therefore $v=y_n-u$ is a weak solution of $v_t={\hcal}(v)$. 
\enddemo
Taking into account the semigroup property \ref{SGR} and independence of defining function we deduce from Lemma \ref{EQ20} the following
\bl\label{EQ21}
Let $v=v(z,t)$ be a local weak solution of $v_t=\hcal(v)$. Suppose that, locally at $(z^0,t^0)$, $v=0$ is a graph
$y_n=u(x,y',t)$ of a contiunous function. Then $u$ is a weak solution of $u_t={\hcal}_0(u)$.
\el
Now we are in position to prove the
following
\bt\label{EQ22}
Let $D$ be a bounded strictly pseudoconvex domain in $\C^{n-1}\tms\R$, $\Gamma_0$ the graph of a continuous function $u_0:{\overline
D}\to\R$. Then the evolution of
$\Gamma_0$
with fixed boundary is governed by the following parabolic problem
\begin{equation}\label{EQ40}
\begin{cases}
u_t={\hcal}_0(u) \>\> & in\>\> D\times (0,+\IN)\\
u(x,y',0)=u_0(x,y')\>\> & {\rm for}\>\>(x,y',0)\in{\overline D}\tms \{0\}\\
u(x,y',t)=u_0(x,y')\>\> & {\rm for}\>\>(x,y',t)\in{\rm b}\,D\tms [0,+\IN).
\end{cases}
\end{equation}
\et
\demo
Let the evolution be defined by the zero set $\{v=0\}$ where $v$ is the weak solution of the parabolic problem $(P)$.

In view of Theorem \ref{EQ25} every $E_t(\Gamma_0,{\rm b}\,\Gamma_0)$, $t\ge 0$, is a graph, a priori over a subset of $\overline D$, but in view of Theorem \ref {EQ16} it separates $D\tms i\R$ so is the graph over $\overline D$, say of a continuous function $u^t=u^t(x,y')$. Define
$u:\overline D\tms(0,+\IN)\to\R$ by $u(x,y',t)=u^t(x,y')$. The function $u$ is continuous: if $(x^n,y'^n,t^n)\to
(\bar x,\bar y',\bar t)$ then the sequence $\left(x^n,y'^n,u^{t_n}(x^n,y'^n),t^n\right)$ tends to a point
$(\bar x,\bar y',\bar y_n,\bar t)$ which lies on the graph of $u^t$. In particular
$\bar y_n=u(\bar x,\bar y',\bar t)$.

Thus
$$
E_t(\Gamma_0,{\rm b}\,\Gamma_0)=\{y_n=u(x, y',t)\}.
$$
Owing to Lemma \ref{EQ21} $u$ is a weak solution of $u_t={\hcal}_0(u)$ which
satisfies all conditions (\ref {EQ40}). This concludes the proof. 
\enddemo

The following lemma will be used in the next section
\bl\label{EQ55}
Let $U$ be a domain in $\C^n$ and $u\in C^0\left(U\tms(0,+\infty)\right)$ a continuous subsolution of $u_t=\hcal_0(u)$ such that $\hcal_0(u)\le 0$ (in the weak sense). Then $u$ is non increasing in time. 
\el
\demo
This follows from the more general fact: let $W=V\times (a,b)\subset \R^N$, $V$ open in $\R^{N-1}$, $u=u(x,t)$ an upper semicontinuous function in $W$ such the inequality $u_t\le 0$ is satisfied in $W$ (in the weak sense). Then, for every $x\in V$,
$a<t_1<t_2<b$ we have $u(x,t_1)\ge u(x,t_2)$. 

Fix $t_1$ and let
$$ W_1=\left\{(x,t)\in\R^N:x\in V, t_1<t<b\right\}.
$$ We may assume, without loss of generality that 
$u\le M<+\IN$ on $W_1$ ($M$ constant), $m=\inf\;g>-\IN$ and that $b-t_1<1$. It
suffices to show the following: for every $v\in C^\IN (V)$ such that 
$u(x,t_1)<v(x)$, it holds
$u(x,t)<v(x)$ for every $(x,t)\in W_1$. 

Set, for $\a\in [0,+\IN)$,
$$
\varphi^\a(x,t)=v(x)+(M-m)(t-t_1)^\a.
$$
Then $\varphi^\a\in C^\IN (W_1)$ and 
$$
\varphi^0 (x,t)=v(x)+(M-m)\ge u(x,t),
$$
$$
\lim_{\a\rightarrow +\IN}\varphi^\a(x,t)=v(x),
$$ 
for $t-t_1<1$, $(x,t)\in W_1$.

Suppose now that $u(x_0,t_0)>v(x_0)$ for some
$(x_0,t_0)\in W_1$. Then there is an $\a\in (0,+\IN)$ and $(x^\ast,t^\ast)\in
W_1$ such that
$\varphi^\a(x^\ast,t^\ast)=u(x^\ast,t^\ast)$. Since $u_t\le 0$ in the weak sense, 
$\varphi^\a((x^\ast,t^\ast)\le 0$, a contradiction.  
\enddemo
\section{Limit for solutions}
In order to describe the asymptotic behaviour of the weak solution $u$ of (\ref{EQ40}) we need to recall some results about the existence of Levi flat hypersurfaces with prescribed boundary.

Let $S\sbs\C^n$ be a connected smooth submanifold of dimension $(2n-2)$. Assume that:
\begin{enumerate}
\item $S$ is compact and nowhere minimal at its {\rm CR} points;
\item $S$ has at least one complex point and every such point of  is flat and elliptic;
\item $S$ does not contain complex manifold of dimension $(n-2)$.
\end{enumerate}
Then in \cite{DTZ1} the following two theorems are proved
\bt\label{kon}
$S$ is diffeomorphic to the unit sphere with two complex points $p_1, p_2$. The CR orbits of $S$ are topological $(2n-3)$-spheres that can be represented as level sets of a smooth function $\nu:S\to\R$, inducing on $S_0=S\setminus\{p_1, p_2\} $ a foliation $\mathcal F$ of class $C^\infty$ with $1$-codimensional compact leaves. 
\et
\bt\label{kon1}
There exist a smooth submanifold $\widetilde S$ and a Levi-flat $(2n-1)$-subvariety $\widetilde M$ in $\R\times\C^n$ {\rm (}i.e. $\widetilde M$ is Levi-flat in $\C\times\C^n${\rm)}, both contained in $[0,1]\times\C^n$, such that $\widetilde S=d{\widetilde M}$ in the sense of currents and the natural projection $\pi\colon [0,1]\times\C^n\to \C^n$ restricts to a diffeomorphism between $\widetilde S$ and $S$.
\et
We can go further if $S$ is a graph:
\bt{\rm (\cite[Theorem 3.1]{DTZ})}
Let $D\sbs\C^{n-1}\tms\R$ be a strictly pseudoconvex bounded domain, $g_0:{\rm b}\,D\to\R$ a smooth function. Assume that the graph $S$ of $g_0$ satisfies the hypothesis of Theorem \ref{kon1}. Then there exists a Lipschitz function $f:\oli D\rightarrow \R$ which is smooth on $\oli D\setminus \{q_1,q_2\}$, the projections of the only two complex elliptic points of $S$, and such that $f_{\vert {\rm b}\,D}=g_0$ and $M={\rm graph}(f)\smi S$ is a Levi flat hypersurface of $\C^n$. 
\et
We want to prove that in this situation, the evolution of an arbitrary, continuous graph over $D$, with boundary $S$ tends as $t\to+\infty$ to the Levi flat graph $M$. This follows from
\bt\label{EQ50}
Let $u_0:\oli D\to\R$ be a continuous function such that $g_0={u_0}_{\vert{\rm b}\,D}$ and $u\in C^0\left({\overline D}\times [0,+\IN)\right)$ the weak solution of the problem {\rm (}\ref{EQ40}{\rm )}. Then
$$
\lim\limits_{t\to+\IN}u(\cdot,t)=f
$$
in $ C^0(\overline D)$. In particular, if $\G_0={\rm graph}{\rm (}u_0{\rm )}$ we have
$$
E_t(\G_0,S)\to M
$$ as
$t\to+\IN$ in the $ C^0$-topology.
\et
\demo
Observe that $\hcal_0(f)=0$ on $D$, $M={\rm graph}(f)\smi S$ being a Levi flat hypersurface.
We divide the proof in several steps. First of all we construct two smooth barriers $\d^\pm$: $\d^-\le u_0\le\d^+$ in $D$, $\d^-=\d^+=u_0$ on ${\rm b}\,D$ and ${\hcal}_0(\d^-)\ge 0$, $\hcal_0(\d^+)\le 0$ in $D$. This is easily done using the functions
$\d^\pm=u_0\mp\lambda\r$ where $\lambda$ is a suitable positive constant $\r=\r(x,y')$ a smooth function on a neighbourhood $U$  of $\overline D$ with the following properties: $\r$ is strictly plurisubharmonic in $D\tms i\R$, $D=\{\r<0\}$ and $d\r\neq 0$ on ${\rm b}\,D$.  

Next we consider the weak solutions $u^\pm$ of (\ref{EQ40}) corresponding respectively to the boundary values $\d^\pm$ on $D$ and $g_0$ on ${\rm b}\,D\tms[0,+\IN)$. $u^\pm$ are bounded by virtue of the maximum principle and uniformly continuous because
of Walsh's Lemma. Moreover, by Lemma \ref{EQ55}, $u^+$ ($u^-$) is non increasing (non decreasing) in $t$
since $\hcal_0(\d^+)\le 0$  (${\hcal}_0(\d^-)\ge 0$). It follows that
$\lim\limits_{t\to+\IN}u^\pm(\xi,t):=\tilde{u}^\pm(\xi)$ exists pointwise.

Now define
functions $u^\pm_h(\cdot,t)=u^\pm(\cdot,t+h)$ for each positive $h$. These
functions are
still weak solutions (with different boundary values). Moreover, since
$u^\pm$ are
bounded, the sets $\left\{u^\pm_h\right\}_{h\ge 0}$ are equicontinuous and  $$
\tilde{u}^\pm(x)=\lim_{t\to+\IN}u^\pm(x,t)=\lim_{t\to+\IN}u^\pm_h(x,t) $$
for every
$\xi\in D$. It follows that $\tilde{u}^\pm$ are continuous in $D$,
$\tilde{u}^\pm=\phi_0$ on ${\rm b}D$ and ${\hcal}_0(u^\pm)=0$ in $D$ and
consequently
(by uniqueness) $\tilde {u}^+=\tilde {u}^-=w$ in $D$.
Consider now the weak solution $u$
of the parabolic problem (\ref{EQ40}). By virtue of the comparison principle we
have 
$$
u^-(\cdot,t)\le u(\cdot,t)\le u^+(\cdot,t)
$$ 
and from this, letting
$t\to+\IN$ we
obtain
$$
f(\xi)=\lim_{t\to+\IN}u^-(\xi,t)\le \liminf_{t\to+\IN}u(\xi,t)\le\limsup_{t\to+\IN}
u(\xi,t)=\lim_{t\to+\IN}u^+(\xi,t)=f(\xi)
$$
for every $\xi\in D$, so
$$
\liminf_{t\to+\IN}u(\cdot,t)=\limsup_{t\to+\IN}u(\cdot,t)=\lim_{t\to+\IN}u(\cdot
,t)=f
$$
in $C^0(\oli D)$.  
\enddemo

\end{document}